\documentclass[reqno,a4paper,11pt]{amsart}

\numberwithin{equation}{section}
\usepackage[utf8]{inputenc}
\usepackage{amsmath, amsthm, amssymb, amscd, accents, bm, amsfonts}
\usepackage{mathtools}
\usepackage{mathrsfs,dsfont}
\usepackage{url}
\usepackage{mathtools}
\usepackage{mdwlist}
\usepackage{paralist}
\usepackage{datetime}
\usepackage{hyperref}
\usepackage{colonequals}
\usepackage[sort,nocompress]{cite}
\mathtoolsset{showonlyrefs}

\newtheorem{definition}{Definition}[section]
\newtheorem{theorem}[definition]{Theorem}
\newtheorem{lemma}[definition]{Lemma}
\newtheorem{corollary}[definition]{Corollary}
\newtheorem{proposition}[definition]{Proposition}
\newtheorem{problem}[definition]{Problem}
\newtheorem{example}[definition]{Example}
\newtheorem{remark}[definition]{Remark}

\def\N{{\mathbb N}}
\def\Z{{\mathbb Z}}
\def\R{{\mathbb R}}
\def\T{{\mathbb T}}
\def\C{{\mathbb C}}

\newcommand{\E}{{\mathbf{E}}}
\newcommand{\Rd}{{\R^d}}
\newcommand{\Cd}{{\C^d}}
\newcommand{\lspan}{{\mathrm{span} \,}}
\newcommand{\supp}{{\mathrm{supp}}}

\newcommand{\ltd}{{L^2(\R^d)}}

\newcommand{\re}{{\mathrm{Re} \,}}
\newcommand{\ift}{{\mathcal{F}^{-1}}}
\newcommand{\ft}{{\mathcal{F}}}
\newcommand\1{\mathds{1}}

\begin{document}

\title[From completeness of translates to phaseless sampling]{From completeness of discrete translates to phaseless sampling of the short-time \\ Fourier transform}

\author[Philipp Grohs]{Philipp Grohs}
\address{Faculty of Mathematics, University of Vienna, Oskar-Morgenstern-Platz 1, 1090 Vienna, Austria \textbf{and} Research Network DataScience@UniVie, University of Vienna, Kolingasse 14-16, 1090 Vienna, Austria \textbf{and} Johann Radon Institute of Applied and Computational Mathematics, Austrian Academy of Sciences, Altenbergstrasse 69, 4040 Linz, Austria}
\email{philipp.grohs@univie.ac.at}

\author[Lukas Liehr]{Lukas Liehr}
\address{Faculty of Mathematics, University of Vienna, Oskar-Morgenstern-Platz 1, 1090 Vienna, Austria}
\email{lukas.liehr@univie.ac.at}

\author[Irina Shafkulovska]{Irina Shafkulovska}
\address{Faculty of Mathematics, University of Vienna, Oskar-Morgenstern-Platz 1, 1090 Vienna, Austria \textbf{and} Acoustics Research Institute, Austrian Academy of Sciences, Wohllebengasse 12-14, 1040 Vienna, Austria}
\email{irina.shafkulovska@univie.ac.at}

\date{\today}

\subjclass[2020]{41A30, 42C30, 94A12, 94A20}
\keywords{phaseless sampling, completeness of translates, uniqueness problem}

\begin{abstract}
We study the uniqueness problem in short-time Fourier transform phase retrieval by exploring a connection to the completeness problem of discrete translates. Specifically, we prove that functions in $L^2(K)$ with $K \subseteq \Rd$ compact, are uniquely determined by phaseless lattice-samples of its short-time Fourier transform with window function $g$, provided that specific density properties of translates of $g$ are met. 
By proving completeness statements for systems of discrete translates in Banach function spaces on compact sets, we obtain new uniqueness statements for phaseless sampling on lattices beyond the known Gaussian window regime. Our results apply to a large class of window functions, which are relevant in time-frequency analysis and applications.
\end{abstract}

\maketitle

\section{Introduction}

The short-time Fourier transform (STFT) of a function $f \in \ltd$ with respect to a window function $g \in \ltd$ is given by
$$
V_gf(x,\xi) = \int_\Rd f(t) \overline{g(t-x)} e^{-2\pi i \xi \cdot  t} \, dt,
$$
with $\xi \cdot t$ denoting the Euclidean inner product in $\Rd$.
For a set $U \subseteq \R^{2d}$, the phaseless STFT samples at $U$ are defined by
$$
|V_gf(U)| \in [0,\infty)^{U}, \quad |V_gf(U)| \coloneqq \left \{ |V_gf(u)| \right  \}_{u \in U}.
$$
The \emph{STFT phase retrieval problem} concerns the inversion of the nonlinear operator
$
f \mapsto |V_gf(U)|
$
that maps a square-integrable function $f$ to its phaseless STFT samples at $U$ with respect to the window function $g$. This problem is significant in various applications, particularly in diffraction imaging \cite{zhou2020low,nature2} and quantum mechanics \cite{orl1994phase}. It has been studied from various perspectives recently \cite{GrohsLiehrSqrt,filbir,ALAIFARI202134,Iwen2022,bartusel2023injectivity,GrohsRathmair,Luef2019}. The question of whether $f \mapsto |V_gf(U)|$ is injective, or, in other words, whether every $f$ can be recovered from $|V_gf(U)|$ in a unique way, is a ubiquitous problem in this research field. This question, known as the \emph{uniqueness problem in STFT phase retrieval}, is especially relevant for computational applications when $U$ is a discrete set. Clearly, if $h=cf$ for some $\T = \{ z \in \C : |z|=1 \}$, then $|V_gf(U)| = |V_gh(U)|$. This implies that uniqueness must be considered modulo a multiplicative constant $c \in \T$. The equivalence relation $f \sim h$ indicates that $h=cf$ for some $c \in \T$. The latter considerations lead to the notion of a uniqueness set.

\begin{definition}\label{def:us}
    A set $U \subseteq \R^{2d}$ is said to be a uniqueness set for phase retrieval in $X \subseteq \ltd$ with window $g \in \ltd$ if the map
    $$
    f \mapsto |V_gf(U)|
    $$
    is injective on $X / \hspace{-0.13cm} \sim$, i.e., if $f,h \in X$ satisfy $|V_gf(u)| = |V_gh(u)|$ for all $u \in U$, then $f = ch$ for some $c \in \T$.
\end{definition}

In view of paradigms in time-frequency analysis and numerical feasibility, it is desirable to select $U$ as a discrete set of sampling locations \cite{Groechenig}. A particularly important and widely studied choice for $U$ are lattices, i.e, $U=M\Z^{2d} = \{ Mz : z \in \Z^{2d} \}$ for some invertible matrix $M \in \mathrm{GL}(2d,\R)$, known as the generating matrix of $U$.

It was recently shown that if $g$ is an arbitrary window function and $X=\ltd$, then every (generic) lattice does not serve as a uniqueness set for phase retrieval in $\ltd$ with window $g$ \cite{grohsLiehr4}. This result implies that the space $\ltd$ is too large to ensure uniqueness from lattice samples. Moreover, it raises the natural question of whether uniqueness can be achieved for function spaces $X$ that are proper subspaces of $\ltd$.

An initial result addressing this question was presented in \cite{grohsliehr1} which forms the basis for the current study. It was proved that for the interval $K=[-\frac{c}{2},\frac{c}{2}]$ and the Gaussian $g(t)=e^{-\pi t^2}$, it holds that $\frac{1}{b}\Z \times \frac{1}{2c}\Z$ is a uniqueness set for phase retrieval in $L^4[-\frac{c}{2},\frac{c}{2}]$ with window $g$ (for every $b>0$). The proof relies on a theorem by Zalik on the completeness of discrete translates of the Gaussian function \cite{Zalik}. Subsequently, it was shown in \cite{wellershoff2024injectivity} that by adapting the proof techniques from \cite{grohsliehr1}, a similar uniqueness result can be obtained with $L^4[-\frac{c}{2},\frac{c}{2}]$ replaced by $L^2[-\frac{c}{2},\frac{c}{2}]$. These findings provide a three-fold motivation for the present study. Specifically, the goals are to
\begin{enumerate}
    \item[(a)] elaborate on a systematic relation between completeness properties of translates and uniqueness sets for the phase retrieval problem in spaces of compactly supported functions,
    \item[(b)] make the results independent of the dimension, as in various applications one faces the situation where $d \geq 2$,
    \item[(c)] show that uniqueness via sampling on lattices is achievable for a large class of naturally occurring window functions, beyond Gaussians.
\end{enumerate}

The primary focus of the present paper can thus be summarized as the following problem, which forms the core of our investigation.

\begin{problem}\label{problem}
    Let $g \in \ltd$ and let $K \subseteq \Rd$ be compact. Does there exist a lattice $\Lambda \subseteq \R^{2d}$ such that $\Lambda$ is a uniqueness set for phase retrieval in $L^2(K)$ with window $g$?
\end{problem}

We remark that in Problem \ref{problem}, the space $L^2(K)$ is identified with the subspace $\{ f \in \ltd : \supp(f) \subseteq K \} \subseteq \ltd$. In this way, Problem \ref{problem} is consistent with the notion of a uniqueness set for spaces $X \subseteq \ltd$ given in Definition \ref{def:us}.

\section{Main results}

\subsection{Uniqueness-completeness relation}

Our first main result establishes a link between Problem \ref{problem}, completeness problems of exponentials, and completeness problems of discrete translates. It asserts that an interplay between completeness of exponentials and completeness of translates gives rise to uniqueness results for the phase retrieval problem and a positive resolution of Problem \ref{problem}.

To state this precisely, we introduce the following notation. For a function $f : \Rd \to \C$ and $\lambda\in \Rd$, the translate of $f$ by $\lambda$ is the function $T_\lambda f \coloneqq f(\cdot - \lambda)$. Given a set $\Lambda \subseteq \Rd$, the system of $\Lambda$-translates of $f$ is defined as
$$
\mathcal{T}(f,\Lambda) \coloneqq \{ T_\lambda f : \lambda \in \Lambda \}.
$$
If $(X,\|\cdot \|)$ is a Banach space of functions defined on a measurable set $K$, and if $f$ has the property that $f(\cdot - \lambda)|_K \in X$ (the restriction of $f(\cdot - \lambda)$ to $K$ belongs to $X$), then $\mathcal{T}(f,\Lambda)$ is complete in $X$ if its linear $\C$-span is dense in $X$. Completeness properties of translates in function spaces were studied by various authors, see, for instance, \cite{LANDAU1972438,olevskii2,olevskii1}.

We denote by $\mathcal{E}(\Gamma) \coloneqq \{ e^{2\pi i \gamma \cdot t} : \gamma \in \Gamma \}$ the system of complex exponentials on $\Rd$ formed by $\Gamma \subseteq \Rd$. Moreover, the following abbreviation is used throughout the article: for $g : \Rd \to \C$ and $\omega \in \Rd$, the map $g_\omega : \Rd \to \C$ is defined by 
$$
g_\omega(t)=g(t-\omega)\overline{g(t)}.
$$
Finally, for $A,B \subseteq \Rd$, the sets $A+B$ and $A-B$ represent the Minkowski sum and difference, respectively.

\begin{theorem}\label{thm:relation}
    Let $K\subseteq \Rd$ be compact and let $\Gamma \subseteq \Rd$ such that $\mathcal{E}(\Gamma)$ is complete in $L^2(K-K)$. Further, let $g \in C(\Rd)$ and let $\Lambda \subseteq \Rd$ such that $\mathcal{T}(g_\omega,\Lambda)$ is complete in $C(K)$ for every $\omega \in K-K$. Then $\Lambda \times \Gamma$ is a uniqueness set for phase retrieval in $L^2(K)$ with window $g$.
\end{theorem}

One can show (see Lemma \ref{lma:diam}) that the choice $\Gamma = \frac{1}{\mathrm{diam}(K)}\Z^d$ implies that $\mathcal{E}(\Gamma)$ is complete in $L^2(K-K)$. Here, $\mathrm{diam}(K)$ denotes the diameter of $K$,
$$
\mathrm{diam}(K) \coloneqq \sup \{ |k-k'|: k,k' \in K \},
$$
with $|\cdot |$ the Euclidean norm in $\Rd$.
It therefore follows from Theorem \ref{thm:relation} that if $\Lambda=A\Z^d$ is a lattice and $g$ is a window function such that $\mathcal{T}(g_\omega,\Lambda)$ is complete in $C(K)$ for every $\omega \in K-K$, then the $2d$-dimensional lattice
$$
U = A\Z^d \times \tfrac{1}{\mathrm{diam}(K)}\Z^d
$$
is a uniqueness set for phase retrieval in $L^2(K)$ with window $g$. 
Consequently, Problem \ref{problem} is effectively addressed through the lens of a completeness property of translates of the window function $g$. Our focus lies therefore on the establishment of such completeness results with respect to translation sets of the form $\Lambda=A\Z^d$ (we remark, however, that most results extend to more general types of sets).

\subsection{Gaussians}\label{sec:exp_decay}

A theorem of Zalik implies that if $\varphi(t) = e^{-ct^2}$ is a Gaussian with $c > 0$, then for any interval $K \subseteq \mathbb{R}$ and any $a>0$, the system $\mathcal{T}(\varphi, a\Z)$ is complete in $L^2(K)$ \cite[Thm. 4]{Zalik}. It was shown in \cite{wellershoff2024injectivity} that this statement also holds for $L^2(K)$ replaced by $C(K)$. These statements were used to establish phaseless sampling results, as detailed in \cite{grohsliehr1} and extended in \cite{wellershoff2024injectivity}. The first consequence of Theorem \ref{thm:relation} yields a multivariate version of the results in \cite{grohsliehr1,wellershoff2024injectivity}.

\begin{proposition}\label{prop:gauss}
Let $g \in C(\Rd)$ be the Gaussian $g(x) = e^{-x^TAx}$ with $A \in \C^{d \times d}$ such that $\re A$ is positive definite. Further, let $K \subseteq \Rd$ be compact. If $\Lambda \subseteq \Rd$ is an arbitrary lattice, then
$$
U = \Lambda \times \tfrac{1}{\mathrm{diam}(K)}\Z^d
$$
is a uniqueness set for phase retrieval in $L^2(K)$ with window $g$.
\end{proposition}

The set $U$ in Proposition \ref{prop:gauss} depends on the set $K$.
It will be pointed out in Corollary \ref{cor:kahane}, that upon replacing the set $\tfrac{1}{\mathrm{diam}(K)}\Z^d$ by a suitable irregular set of sampling locations $\Gamma \subseteq \Rd$, the uniqueness set in Proposition \ref{prop:gauss} can be made independent of $K$. That is, one can find an irregular set $\Gamma$ such that $\Lambda \times \Gamma$ is a uniqueness set for the space
$$
X = \bigcup_{\substack{K \subseteq \Rd \\ K \, \mathrm{compact}}} L^2(K).
$$
Notice, that $X$ is dense in $\ltd$. We further emphasize that $\Gamma$ can be chosen to very sparse, in the sense that the point-density of $\Gamma$ is equal to zero. We refer to Section \ref{sec:gauss} for further discussions.

\subsection{Bandlimited functions}

A function $g \in \ltd$ is said to be bandlimited to the compact set $K' \subseteq \Rd$ if the support of the Fourier transform of $g$ is contained in $K'$. The Paley-Wiener space $PW_{K'}$ is the collection of all functions in $\ltd$ that are bandlimited to $K'$. It is known that suitable lattice-translates of a bandlimited function are complete in $C(K)$ (this is well-known for the univariate case $d=1$, see \cite[Prop. 5.5]{pinkus}; the statement readily extends to the multivariate case, as shown in the proof of Proposition \ref{prop:bandlimited}). Combining this result with Theorem \ref{thm:relation} provides a uniqueness statement for phaseless sampling on lattices $U=M\Z^{2d}$ with bandlimited window functions. Notably, in this regime, the lattice $U$ does not need to be separable (i.e. $U=A\Z^d \times B \Z^d$ for some $A,B \in \mathrm{GL}(d,\R)$) but merely a condition on the dual lattice $U^* = M^{-T}\Z^{2d}$ must be satisfied.

\begin{proposition}\label{prop:bandlimited}
    Let $K,K' \subseteq \Rd$ be compact sets, and let $g \in PW_{K'}$, $g \neq 0$. If $U \subseteq \R^{2d}$ is a lattice such that $(K'-K') \times (K-K)$ is contained in a fundamental domain of the dual lattice $U^*$, then $U$ is a uniqueness set for phase retrieval in $L^2(K)$ with window $g$. For instance, this holds if
    $$
    U = \tfrac{1}{\mathrm{diam}(K')}\Z^d \times \tfrac{1}{\mathrm{diam}(K)}\Z^d.
    $$
\end{proposition}

We remark, that if $g \in PW_{K'}$ and if $f \in L^2(K)$, then $V_gf \in PW_{K' \times (-K)}$. Hence, the image of $L^2(K)$ under $V_g$ is a subspace of $PW_{K' \times (-K)}$. It follows that Proposition \ref{prop:bandlimited} admits a reformulation as a uniqueness statement in infinite-dimensional subspaces of Paley-Wiener spaces of the form $PW_S$ where $S$ is a separable compact subset of $\R^{2d}$ (i.e. $S=A\times B$ for compact subsets $A,B \subseteq \Rd$). The uniqueness-theory of phase retrieval in (subspaces) of Paley-Wiener spaces was previously studied in various works \cite{lai2021conjugate,zhang2024gabor,aku1,Jaming2022}. We refer to Remark \ref{rem:pw} for a detailed comparison of Proposition \ref{prop:bandlimited} to the literature on phase retrieval in Paley-Wiener spaces.

\subsection{Hermite functions and the class $\mathcal{P}_{\alpha,\beta}$}
In various contexts of time-frequency analysis and applications, one deals with window functions that are Hermite functions, see, for instance \cite{horst2023non,grochenig2009gabor,grochenig2007gabor} and the references therein. In what follows, we consider a more general setting that covers these windows. To do so, we recall that an entire function $h : \Cd \to \C$ is said to be of exponential type smaller than or equal to $\alpha$ if there exists a constant $C>0$ such that 
$$
|h(z)| \leq C e^{\alpha \| z \|_1}, \quad z \in \Cd,
$$
where $\| z \|_1 = \sum_{j=1}^d |z_j|$ for $z = (z_1, \dots, z_d) \in \C^d$.
Let $\E_\alpha(\Cd)$ denote the collection of all entire functions of exponential type smaller than or equal to $\alpha$. Then we define the function class $\mathcal{P}_{\alpha,\beta}$ by
$$
\mathcal{P}_{\alpha,\beta} \coloneqq \left \{ p(z)e^{-z^TAz+b\cdot z} : p \in \E_\alpha(\Cd), \, A \in \C^{d \times d}, \, \|A + A^T \|_1 \leq \beta, \, b\in\Cd \right \}.
$$
Here, $\| M \|_1 = \max_{j=1,\dots,d} \sum_{i=1}^d |a_{ij}|$ denotes the maximal absolute column sum of a matrix $M \in \C^{d \times d}$.

Since every polynomial is an element of $\E_\alpha(\Cd)$ for every $\alpha >0$, and every Hermite function is a product of a Gaussian with a polynomial, it follows that $\mathcal{P}_{\alpha,\beta}$ contains the Hermite functions for suitable $\alpha,\beta>0$. In regards of Problem \ref{problem} with windows belonging to $\mathcal{P}_{\alpha,\beta}$, the following can be said.

\begin{proposition}\label{prop:p_ab}
    Let $K \subseteq \R^d$ be compact, and let $g \in \mathcal{P}_{\alpha,\beta} \cap \ltd$, $g \neq 0$. If $c \in (0,\infty)$ satisfies $c > \frac{2\alpha +\mathrm{diam}(K)\beta}{\pi}$, then 
    $$
    U = \tfrac{1}{c}\Z^d \times \tfrac{1}{\mathrm{diam}(K)}\Z^d
    $$
    is a uniqueness set for phase retrieval in $L^2(K)$ with window $g$.
\end{proposition}

\section{Proof and consequences of Theorem \ref{thm:relation}}

Before we begin to prove the main results, we introduce the following notation.

The Fourier transform $\ft f$ of $f\in \ltd$ is the unitary operator given by
\begin{equation}\label{eq:def_ft}
    \ft f(\omega) = \int_{\Rd} f(t) e^{-2\pi i \omega\cdot t}dt,\qquad \omega\in\Rd.
\end{equation}
The integral in \eqref{eq:def_ft} is well-defined whenever $f\in L^1(\Rd)\cap \ltd$ and the operator extends to $\ltd$ in the usual way.
Furthermore, applying the Fourier transform twice yields the reflection operator: $\ft^2 f= \mathcal{R} f \coloneqq f(-\cdot)$. Finally, we mention the identities
\begin{equation}\label{eq:fourier_identities}
    \ft f = \ift \mathcal{R}f = \mathcal{R} \ift f, \quad \overline{\ft f} = \ift \overline{f}, \quad \ft \overline{f} = \overline{\mathcal{R}\ft f}.
\end{equation}
The following elementary lemma will be used frequently in the upcoming proofs.

\begin{lemma}\label{lem:supp_of_squares}
    Let $f \in \ltd$ and let $K = \supp(\ft f)$. Then $\supp(\ft |f|^2 ) \subseteq K-K$.
\end{lemma}
\begin{proof}
    By the convolution theorem and \eqref{eq:fourier_identities},
    \begin{equation}\label{eq:supp_lemma}
        \begin{split}
            \ft |f|^2 = \ft({f}\cdot \bar{f}) = \ft f * \ft \bar f =  \ft f * \overline{\mathcal{R}\ft f}. 
        \end{split}
    \end{equation}
    By assumption, $\overline{\mathcal{R}\ft f}$ is supported on $-K$. Thus, the support of the convolution in \eqref{eq:supp_lemma} is contained in $K + (-K) = K-K$.
\end{proof}

We continue with the proof of the first main result of the paper.

\begin{proof}[Proof of Theorem \ref{thm:relation}]
    Let $f \in L^{2}(K)$. Since $g \in C(\Rd)$, it follows that $f\overline{T_xg} \in L^2(K)$ for every $x \in \Rd$. The properties of the Fourier transform in \eqref{eq:fourier_identities} and the definition of the STFT yield
    \begin{equation*}
    \begin{split}
        |V_gf(x,\cdot)|^2 & = \ft(f\,\overline{T_xg}) \, \overline{\ft(f\,\overline{T_xg})} = \ift \mathcal{R} (f\,\overline{T_xg})\,\ift (\overline{f}\,T_xg) \\
        & = \ift \left( \mathcal{R} \left(f\,\overline{T_xg}\right) * \left(\overline{f}\, T_xg\right)\right).
    \end{split}
\end{equation*}    
As $f$ is supported on $K$, we obtain the identity
\begin{equation*}
    \begin{split}
        \mathcal{R} \left(f\,\overline{T_xg}\right) * \left(\overline{f}\,T_xg\right)(\omega) & = \int_K f(-(\omega-t)) \,\overline{g(-(\omega-t)-x)}\, \overline{f(t)}\, g(t-x) \, dt \\
        & = \int_K f_\omega(t) \overline{T_x g_\omega(t)} \, dt.
    \end{split}
\end{equation*}
Therefore, it holds that
\begin{equation}\label{eq:ft_identity}
    \ft \left(|V_gf(x,\cdot)|^2\right)(\omega) = \int_K f_\omega(t) \overline{T_x g_\omega(t)} \, dt.
\end{equation}

Since the support of $f$ is contained in $K$, and $V_gf(x,\cdot) = \ft (fT_x \overline{g})$, it follows that $\ft \left(V_gf(x,\cdot)\right)$ has support in $-K$. According to Lemma \ref{lem:supp_of_squares}, $\ft \left(|V_gf(x,\cdot)|^2\right)$ has support in $K-K$.

Further, since $f\,\overline{T_xg} \in L^{2}(K)$, we have
\begin{equation}
    |\ft (f\overline{T_xg})|^2 = |V_gf(x,\cdot)|^2 \in L^1(\Rd).
\end{equation}
Hence, $\ft ( |V_gf(x,\cdot)|^2)$ is a continuous function supported in $K-K$ (this also yields $\ft ( |V_gf(x,\cdot)|^2) \in L^2(K-K)$). The completeness of $\mathcal{E}(\Gamma)$ in $L^2(K-K)$ implies that $|V_gf(x,\cdot)|^2$ is determined by its values on $\Gamma$.

Now let $h \in L^{2}(K)$ be another function whose phaseless STFT samples agree with those of $f$ on the set $\Lambda \times \Gamma$, i.e.,
\begin{equation}\label{eq:lg_eq}
    |V_gf(\Lambda \times \Gamma)| = |V_gh(\Lambda \times \Gamma)|.
\end{equation}
We want to show that there exists $c \in \T$ such that $f=ch$.

As discussed above, $|V_gf(x,\cdot)|^2$ and $|V_gh(x,\cdot)|^2$ are determined by their values on $\Gamma$. Hence, the relation \eqref{eq:lg_eq} implies that
\begin{equation}
    |V_gf(\Lambda \times \Rd)|^2 = |V_gh(\Lambda \times \Rd)|^2.
\end{equation}
An application of the Fourier transform in combination with equation \eqref{eq:ft_identity} shows that
\begin{equation}\label{eq:conv_relation}
    \int_K f_\omega(t) \overline{T_\lambda g_\omega(t)} \, dt = \int_K h_\omega(t) \overline{T_\lambda g_\omega(t)} \, dt \quad \lambda \in \Lambda, \ \omega\in \R^d.
\end{equation}
Now consider 
\eqref{eq:conv_relation}
with $\omega \in \Rd$ fixed. By the Cauchy-Schwarz inequality, it holds that $f_\omega, h_\omega \in L^1(K)$. Hence, the system of translates $\mathcal{T}(g_\omega,\Lambda)$ annihilates the complex regular Borel measure $\mu$ given by $\mu \coloneqq (\overline{f_\omega} - \overline{h_\omega})m$ with $m$ the Lebesgue measure on $K$. 
By assumption, $\mathcal{T}(g_\omega,\Lambda)$ is complete in $C(K)$. The Riesz representation theorem \cite[Thm. 6.19]{rudin} implies that $\mu = 0$, thus $f_\omega = h_\omega$. Since $\omega$ was arbitrary, we can apply the Fourier transform and obtain $\ft(f_\omega)(x) = \ft(h_\omega)(x)$ for every $x \in \Rd$ and every $\omega \in \Rd$. Written out with respect to the STFT, the latter means that $V_ff = V_hh$. It is known that this yields $f=ch$ for some $c \in \T$ (see \cite[Theorem 2.5]{auslander}).
\end{proof}

In order to address Problem \ref{problem}, the following lemma is of importance.

\begin{lemma}\label{lma:diam}
    If $K \subseteq \Rd$ is a compact set and $\Gamma = \frac{1}{\mathrm{diam}(K)}\Z^d$, then $\mathcal{E}(\Gamma)$ is complete in $L^2(K-K)$.
\end{lemma}
\begin{proof}
    Since $K$ is compact, there exists a $k_0$ such that $K_0\coloneqq K-k_0$ is centered around zero, i.e.,
    \begin{equation}
        K_0\subseteq\left[-\tfrac{\mathrm{diam}(K)}{2}, \tfrac{\mathrm{diam}(K)}{2}\right]^d.
    \end{equation}
    By construction, $K_0-K_0 = K-K$, and combined with the centering of $K_0$, one obtains
    \begin{equation}
        K-K=K_0-K_0\subseteq \left[-\mathrm{diam}(K), \mathrm{diam}(K)\right]^d.
    \end{equation}
    It is known, that the system $\mathcal{E}(\Gamma)$ with $\Gamma = \frac{1}{\mathrm{diam}(K)}\Z^d$ is an orthonormal basis for $L^2(\left[-\mathrm{diam}(K), \mathrm{diam}(K)\right]^d)$. In particular, $\mathcal{E}(\Gamma)$ is complete in $L^2(K-K)\subseteq L^2(\left[-\mathrm{diam}(K), \mathrm{diam}(K)\right]^d)$.
\end{proof}

A combination of the previous lemma with Theorem \ref{thm:relation}, implies that the phase retrieval problem in $L^2(K)$ with window $g$ is discretizable whenever lattice-translates of $g_\omega$ are complete in $C(K)$. To be precise, the following statement holds true.

\begin{corollary}\label{cor:relation}
    Let $K\subseteq \Rd$ be a compact set, let $g \in C(\Rd)$, and let $\Lambda \subseteq \Rd$ be a lattice such that $\mathcal{T}(g_\omega,\Lambda)$ is complete in $C(K)$ for every $\omega \in K-K$. Then the phase retrieval problem in $L^2(K)$ with window $g$ admits a lattice-uniqueness set. Specifically, the lattice $\Lambda \times \tfrac{1}{\mathrm{diam}(K)}\Z^d$ is a uniqueness set for phase retrieval in $L^2(K)$ with window $g$.
\end{corollary}
\begin{proof}
    This is a direct consequence of Theorem \ref{thm:relation} and Lemma \ref{lma:diam}.
\end{proof}

\section{Proof and discussion of Proposition \ref{prop:gauss} -- \ref{prop:p_ab}}

\subsection{Gaussian window functions}\label{sec:gauss}

We start by proving that lattices are uniqueness sets for phase retrieval in $L^2(K)$ with Gaussian windows in any dimension.

\begin{proof}[Proof of Proposition \ref{prop:gauss}]
Let $\omega \in K-K$. Following Theorem \ref{thm:relation}, it suffices to show that the system $\mathcal{T}( g_\omega, \Lambda)$ is complete in $C(K)$ whenever 
$g(x) = e^{-x\cdot Ax}$ for some $A\in\C^{d\times d}$ with positive definite real part, and $\Lambda\subseteq\Rd$ is an arbitrary lattice. Using classical product formulas for Gaussians, it follows that the map $g_\omega$ is given by
\begin{equation}
        g_\omega(x) = C \cdot \exp\left(-2\,(x-\tfrac{\omega}{2})^T(\re A)\, (x-\tfrac{\omega}{2})\right),
\end{equation}
where $C \neq 0$ is a constant depending only on $A$ and $\omega$. Hence, it suffices to show that for every function of the form
$$
h: \Rd \to \C, \quad h(x) = \exp\left(-(x-\nu)^TM(x-\nu)\right)
$$
with $\nu \in \R^d$ and $M \in \R^{d \times d}$ positive definite, it holds that $\mathcal{T}(h,\Lambda)$ is complete in $C(K)$ for every lattice $\Lambda$.

To do so, we expand the term $(x-\nu-\lambda)^TM(x-\nu-\lambda)$ and observe that
$$
T_\lambda h(x) = c(\lambda) a(x) B_\lambda(x),
$$
where
\begin{equation*}
    \begin{split}
        c(\lambda) & = \exp \left(-(\lambda+\nu)^T M (\lambda+\nu)\right),\\
        a(x) & = \exp \left ( - x^TMx + 2\left( M \, x\right)^T \nu \right ), \\
        B_\lambda(x) & = \exp\left(2\left(M\, x\right)^T\lambda\right).
    \end{split} 
\end{equation*}
Notice that the completeness of $\mathcal{T}(h,\Lambda)$ is unaffected if each $T_\lambda h$ is multiplied by a non-zero constant. Since $c(\lambda)\neq 0$, $\mathcal{T}(h,\Lambda)$ is complete in $C(K)$ if and only if $\{ a B_\lambda : \lambda \in \Lambda \}$ is complete in $C(K)$. Moreover, $a$ is a smooth, non-vanishing weighting factor independent of $\lambda$. Hence, $\{ a B_\lambda : \lambda \in \Lambda \}$ is complete in $C(K)$ if and only if $\{ B_\lambda : \lambda \in \Lambda \}$ is complete in $C(K)$. Now let $S$ be the complex linear span of the functions $B_\lambda$,
$$
S \coloneqq \lspan \{ B_\lambda : \lambda \in \Lambda \}.
$$
We show that $S$ satisfies the assumptions of the Stone-Weierstrass theorem. Since $B_\lambda$ is real-valued, the set $S$ is invariant under complex conjugation. Further, for $\lambda, \lambda' \in \Lambda$ holds $B_\lambda(x)B_{\lambda'}(x) = B_{\lambda + \lambda'}(x)$, which is an element of $S$ since $\Lambda$ is a lattice (hence, a group). It remains to show that $S$ separates points. To do so, let $x \neq y \in K$ and $\lambda \in \Lambda$. Then the equality $B_\lambda(x) = B_\lambda(y)$ holds if and only if
$$
\left(M(x-y)\right)^T \lambda = 0.
$$
If this equality holds for all $\lambda\in\Lambda$, then it holds for all $\lambda\in\R^d$ (since $\Lambda$ is a lattice and therefore a spanning set for $\Rd$), implying that 
$$
M(x-y) = 0.
$$
Since $M$ is positive definite and therefore invertible, we have $x = y$, contradicting the assumption that $x \neq y$. In conclusion, the Stone-Weierstrass theorem \cite[Thm.~4.51]{FollandRealAnalysis} implies that the closure of $S$ with respect to the maximum norm is either equal to $C(K)$ or equal to $\{ f\in C(K) : f(x_0) = 0 \}$ for a unique $x_0\in K$. The latter, however, is excluded by the positivity of $B_\lambda$.
\end{proof}

Observe that the uniqueness set $U = \Lambda\times \frac{1}{\mathrm{diam}(K)}\Z^d$ in Proposition \ref{prop:gauss} can have an arbitrarily small density since there is no restriction on the density of $\Lambda$. However, there is a dependence on the compact set $K$ in terms of its diameter. One can circumvent this dependence by introducing irregular sampling. To this end, we recall the following result of Kahane \cite{kahane}.

\begin{theorem}[Kahane]
    There exists a symmetric sequence $\mathcal{K} = \{ \kappa_n : n \in \Z \} \subseteq \R$, $\mathcal{K}=-\mathcal{K}$, with the property that the system of exponentials $\mathcal{E}(\mathcal{K})$ is complete in $C(I)$ for every compact interval $I \subseteq \R$. Moreover, $\mathcal{K}$ satisfies the asymptotic
\begin{equation}\label{eq:kahane_zero}
    \lim_{n \to \infty} \frac{n}{\kappa_n} = 0.
\end{equation}
\end{theorem}

In the following, we refer to $\mathcal{K}$ as \emph{Kahane's sequence} and set $\mathcal{K}^d = \mathcal{K} \times \cdots \times \mathcal{K}$ to be the $d$-fold product of $\mathcal{K}$. Using Kahane's sequence, we obtain the following statement.

\begin{corollary}\label{cor:kahane}
    Let $g \in C(\Rd)$ be the Gaussian $g(x) = e^{-x^TAx}$, $A \in \C^{d \times d}$, where $\re A$ is positive definite. Let $\Lambda \subseteq \Rd$ be an arbitrary lattice and let $\mathcal{K}^d \subseteq \Rd$ be Kahane's sequence in $\Rd$. Then $U \coloneqq \Lambda \times \mathcal{K}^d$ is a uniqueness set for phase retrieval in
    $$
    X = \bigcup_{\substack{K \subseteq \Rd \\ K \, \mathrm{compact}}} L^2(K)
    $$
    with window $g$.
\end{corollary}
\begin{proof}
Let $f,h \in X$ such that $\supp(f)=K_1$ and $\supp(h)=K_2$. Further, let $K = [-c,c]^d\subseteq \Rd$ be a cube in $\Rd$ that contains $K_1$ and $K_2$. 

The proof of Proposition \ref{prop:gauss} above shows that $\mathcal{T}(g_\omega,\Lambda)$ is complete in $C(K)$ for all $\omega\in K-K$. Kahane's theorem combined with an elementary tensor argument shows that $\mathcal{E}(\mathcal{K}^d)$ is complete in $L^2(K-K) \subseteq L^2([-2c,2c]^d)$ for every $c>0$. Consequently, $\Lambda \times \mathcal{K}^d$ is a uniqueness set for phase retrieval in $L^2(K)$ with window $g$. Since $f,h \in L^2(K)$, the statement is proved.
\end{proof}

We notice that the space $X$ in Corollary \ref{cor:kahane} is a dense subspace of $\ltd$. In addition, the uniqueness set $U = \Lambda \times \mathcal{K}^d$ has a point density that is equal to zero, i.e.,
$$
\lim_{r \to \infty} \frac{\#(U \cap [-r,r]^{2d})}{r^{2d}} = 0,
$$
where $\#(S)$ denotes the number of elements in a set $S$. The property that the point density of $U$ is zero follows from the relation \eqref{eq:kahane_zero}. This shows how Theorem \ref{thm:relation} can be used to establish uniqueness results for phase retrieval in dense subspaces with respect to sparse uniqueness sets, provided that the sparsity is measured in terms of the point density.

\subsection{Bandlimited window functions}

In order to prove Proposition \ref{prop:bandlimited}, we require the following lemma, which states that the system of all translates of a non-trivial function in $L^1(\Rd) \cap C(\Rd)$ is complete in $C(\Rd)$ with respect to the topology of uniform convergence on compact intervals \cite[Prop.~6.3]{pinkus}.

\begin{lemma}\label{lma:l1}
    If $f \in L^1(\Rd) \cap C(\Rd)$, $f \neq 0$, then $\mathcal{T}(f,\Rd)$ is complete in $C(K)$ for every compact set $K \subseteq \Rd$.
\end{lemma}

In addition, we make use of the following result on sampling in Paley-Wiener spaces \cite[Lem.~14.3]{Higgins}.

\begin{lemma}\label{lma:sampling}
    Let $K \subseteq \Rd$ be a compact set and let $U \subseteq \Rd$ be a lattice such that $K$ is contained in a fundamental domain of the reciprocal lattice $U^*$. If $f,h \in PW_K$ satisfy
    $$
    f(u) = h(u), \quad u \in U,
    $$
    then $f=h$.
\end{lemma}

\begin{proof}[Proof of Proposition \ref{prop:bandlimited}]
Since $f,g \in \ltd$, it follows that $V_gf \in L^2(\R^{2d})$ \cite[Cor. 3.2.2]{Groechenig}. Denoting by $\ft_{2d}$ and $\ft_d$ the Fourier transform on $L^2(\R^{2d})$ and $\ltd$, respectively, it holds that
$$
\ft_{2d}(V_gf)(x,\omega) = e^{2\pi i x \cdot \omega} f(-\omega) \overline{\ft_d g(x)}.
$$
Hence, the assumption on $f$ and $g$ implies that the support of $\ft_{2d}(V_gf)$ is contained in $K' \times (-K)$.
This, in particular, shows that
$$
V_gf \in PW_{K' \times (-K)}.
$$

According to Lemma \ref{lem:supp_of_squares}, it follows that
$$
\supp(\ft_{2d}(|V_gf|^2)) \subseteq (K' - K') \times (K \times K).
$$
In addition, $\ft_{2d}(|V_gf|^2) \in C(\R^{2d})$. This yields $\ft_{2d}(|V_gf|^2) \in L^2(\R^{2d})$ and therefore $|V_gf|^2 \in L^2(\R^{2d})$, where we used that the Fourier transform is an isomorphism. We therefore obtain
$$
|V_gf|^2 \in PW_{(K'-K') \times (K-K)}.
$$
Lemma \ref{lma:sampling} implies that if $U$ is a lattice such that $(K'-K') \times (K-K)$ is contained in a fundamental domain of $U^*$, and if $f,h \in L^2(K)$ satisfy $|V_gf(U)|=|V_gh(U)|$, then $|V_gf(\R^{2d})|=|V_gh(\R^{2d})|$. 

According to Theorem \ref{thm:relation}, it suffices to prove that $\mathcal{T}(g_\omega,\Rd)$ is complete in $C(K)$ for any $\omega \in K-K$. Fix such an $\omega\in K-K$. Since $g \neq 0$ and extends to a holomorphic function, it follows that $g_\omega \neq 0$. 
But $g_\omega \in L^1(\Rd)$. 
Hence, Lemma \ref{lma:l1} implies that $\mathcal{T}(g_\omega,\Rd)$ is complete in $C(K)$ for any compact $K \subseteq \Rd$. This concludes the proof.
\end{proof}

\begin{example}[Airy disk]
Let $D_a \coloneqq \{ x \in \R^2 : |x| \leq a \}$ be the centered disk of radius $a>0$ in $\R^2$, and let $\1_{D_a}$ be the characteristic function of $D_a$. The square of the Fourier transform of $\1_{D_a}$ is called the Airy disk of radius $a$. It holds 
\begin{equation*}
    \begin{split}
        \ft (\1_{D_a})(\omega) &= \int_0^a r \int_0^{2\pi} e^{-2\pi i r | \omega |_2 \cos(\theta)} \, d\theta dr 
        = 2\pi \int_0^a r J_0(2\pi | \omega |_2 r) \, dr.
    \end{split}
\end{equation*}
The above equation follows from a transformation to polar coordinates and the definition of the Bessel functions
$$
J_n(x) = \frac{1}{2\pi} \int_{-\pi}^{\pi} e^{i(x\sin \tau -n\tau)} \, d\tau, \quad n \in \N_0.
$$
The identity
$$
\int_0^a x J_0(x) \, dx = a J_1(a), \quad a>0,
$$
implies that the Airy disk of radius $a$ is then given by the radial function
$$
\mathcal{A}_a(\omega) \coloneqq (\ft (\1_{D_a})(\omega))^2 = \left ( \frac{a J_1(2\pi | \omega |_2 a)}{| \omega |_2} \right )^2.
$$
It is known that the Airy disk behaves approximately like a Gaussian (except for the decay at infinity). The Airy disk commonly appears in diffraction imaging, particularly when an incoming wavefront passes through a circular aperture (a pinhole) before being diffracted by the object of interest. The phase retrieval problem with the Airy disk as a window function naturally arises in this context \cite[Chap.~8.5.2]{bornAndWolf}. 
Proposition \ref{prop:bandlimited} asserts that this type of phase retrieval problem, which is relevant in physical applications, is discretizable. 
Indeed, by the convolution theorem, the Fourier transform of the Airy function $\mathcal{A}_a$ is given $\ft \mathcal{A}_a = \1_{D_a} * \1_{D_a}$, in particular, $\mathcal{A}_a\in PW_{D_{2a}}$.
\end{example}

\begin{remark}[Phase retrieval in Paley-Wiener spaces]\label{rem:pw}
    Let $S \subseteq \R^k$ be a compact set and let $Y \subseteq PW_S$. A set $U \subseteq \R^k$ is said to be a uniqueness set for phase retrieval in $Y$ if the map
    $$
    f \mapsto |f(U)| \coloneqq \{ |f(u)| \}_{u \in U}
    $$
    is injective on $Y / \hspace{-0.13cm} \sim$. It is known that $U=\R^k$ is not a uniqueness set for phase retrieval in $Y = PW_S$. In particular, the phase retrieval problem in $PW_S$ is not discretizable. On the other hand, it is known that sufficiently dense lattices are uniqueness sets for phase retrieval in
    $$
    Y_\R \coloneqq \{ f \in PW_S : f \, \text{real-valued} \, \}
    $$
    of real-valued functions in $PW_S$ \cite{Thakur,alaifariGrohsDaubechies}. As shown in the proof of Proposition \ref{prop:bandlimited}, one has
    $
    V_gf \in PW_{K' \times (-K)}
    $
    whenever $K,K' \subseteq \Rd$ are compact, $f \in L^2(K)$, and $g \in PW_{K'}$. The statement in Proposition \ref{prop:bandlimited} therefore implies that if $k=2d$, $S=K' \times (-K)$, and $U \subseteq \R^{2d}$ is a lattice such that $(K' - K') \times (K \times K)$ is contained in a fundamental domain of $U^*$, then $U$ is a uniqueness set for phase retrieval in
    $$
    Y_{g,K,K'} = \{ V_gf : f \in L^2(K) \} \subseteq PW_S.
    $$
    We have therefore identified an infinite-dimensional subspace of $PW_S$ consisting of complex-valued functions, that admits a uniqueness set which is a lattice. Since $\R^{2d}$ (and therefore $U$) is not a uniqueness for $PW_S$, it follows that $Y_{g,K,K'}$ is necessarily a proper subspace of $PW_S$. For additional studies on phase retrieval in Paley-Wiener spaces, we refer to \cite{lai2021conjugate,zhang2024gabor,aku1,Jaming2022,Thakur}.
\end{remark}

\subsection{The class $\mathcal{P}_{\alpha,\beta}$}
In order to prove Proposition \ref{prop:p_ab}, we recall Carlson's theorem \cite[p. 58, Thm. 3]{levin1996lectures} on uniqueness sets for entire functions of exponential type. 
Though originally stated in the univariate setting, an elementary tensor argument implies the multivariate version as stated next.

\begin{theorem}[Carlson]\label{thm:carlson}
Let $\sigma>0$ and let $c>\frac{\sigma}{\pi}$. If $H \in \E_\sigma(\Cd)$ satisfies $H(\lambda) = 0$ for all $\lambda \in \tfrac{1}{c}\N^d$, then $H$ vanishes identically.
\end{theorem}

We remark, that the set $\tfrac{1}{c}\N^d$ in Carlson's theorem can be replaced by more general (irregular) point sets that satisfy a density condition \cite{rubel1956necessary}.

We now proceed with the proof of Proposition \ref{prop:p_ab}. 

\begin{proof}[Proof of Proposition \ref{prop:p_ab}]
According to Corollary \ref{cor:relation}, we need to prove that the system of translates $\mathcal{T}(g_\omega,\Lambda)$ is complete in $C(K)$ for every $\omega \in K-K$ whenever $\Lambda = \frac{1}{c}\Z^d$ and $c>\frac{\alpha+2\beta\mathrm{diam}(K)}{\pi}$. We split the proof in two steps. In the first step, we show that if $g \in (\mathcal{P}_{\alpha,\beta}\cap \ltd) \setminus \{ 0 \}$, then $g_\omega = (T_\omega g)\overline{g} \in  (\mathcal{P}_{2\alpha,2\beta} \cap L^1(\Rd)) \setminus \{ 0 \}$. In the second step, we relate the problem of completeness of translates
of $g_\omega$ to uniqueness sets in $\mathbf{E}_\sigma(\Cd)$.

\textbf{Step 1.} 
Let $g(x) = e^{-x\cdot Ax+b\cdot x}p(x)\in (\mathcal{P}_{\alpha,\beta}\cap \ltd) \setminus \{ 0 \}$ for some $A\in\C^{d\times d}$ with $\lVert{A+A^T}\rVert_{1}\leq \beta$, $b\in\Cd$, and $p\in\E_\alpha(\Cd)$.
By the Cauchy-Schwarz inequality, it holds that $g_\omega \in L^1(\Rd)$. Since $g_\omega$ is the product of two non-zero entire functions, it follows that $g_\omega$ does not vanish identically.  Furthermore, for all $x\in\Rd$
\begin{equation}
    \begin{split}
        g_\omega(x) 
        &= e^{-(x-\omega)\cdot A(x-\omega) +b\cdot (x-\omega)}p(x-\omega) \overline{e^{-x\cdot Ax+b\cdot x}p(x)} \\
        &= e^{-\omega\cdot A\omega+b\cdot \omega}e^{-x\cdot (A+\overline{A})x+((A+A^T)\omega+b+\overline{b})\cdot x} p(x-\omega) \overline{p(\bar{x})} \\
        &= e^{-x\cdot \tilde A x +\tilde b\cdot x} \tilde p(x),
    \end{split}
\end{equation}
where
$$
\tilde A = A + \overline{A}, \quad \tilde b = (A+A^T)\omega + b + \overline{b}, \quad \tilde p(z) = e^{-\omega\cdot A\omega +b\cdot \omega} p(z-\omega)\overline{p(\overline{z})}.
$$
By the triangle inequality and the assumption on $A$, we have $\lVert \tilde A +\tilde A^T \rVert_1 \leq 2\beta$.
To finish Step 1, it remains to show that $\tilde p$ is an element of $\mathbf{E}_{2\alpha}(\Cd)$. Since $p$ is an entire function, so is $z \mapsto \overline{p(\overline{z})}$. Hence, $\tilde p$ is entire and satisfies
$$
|\tilde p(z)| \leq C^2 e^{\alpha \| z-\omega \|_1} e^{\alpha \| \overline{z} \|_1} \leq C^2 e^{\alpha \| \omega \|_1} e^{2\alpha \| z \|_1},
$$
where we used that $|p(z)| \leq C e^{\alpha \| z\|_1}$ for some $C>0$. This shows that $\tilde p \in \mathbf{E}_{2\alpha}(\Cd)$. Thus, $g_\omega \in  (\mathcal{P}_{2\alpha,2\beta} \cap L^1(\Rd)) \setminus \{ 0 \}$.

\textbf{Step 2.}
Let $\alpha'=2\alpha$ and let $\beta' = 2\beta$.
Further, let $h\in (\mathcal{P}_{\alpha',\beta'}\cap L^1(\Rd)) \setminus \{ 0 \}$, i.e., $h(x) = e^{-x\cdot Ax+b\cdot x}p(x)$ for some $A\in\C^{d\times d}$ with $\lVert{A+A^T}\rVert_1\leq \beta'$, $b\in\Cd$ and $p\in\E_{\alpha'}(\Cd)\setminus\{0\}$. In the following, we show that $\mathcal{T}(h,\frac{1}{c}\Z^d)$ is complete $C(K)$. 

Let $\mu$ be a complex regular Borel measure in $K$. By Riesz representation theorem, it suffices to show that if
\begin{equation}\label{eq:0_funct}
   0=\int_K h(x-\lambda) d\mu(x), \quad \lambda\in \tfrac{1}{c}\Z^d,
\end{equation}
then $\mu = 0$.
Let $k_0\in \Rd$ such that $K_0 \coloneqq K-k_0$ is contained in the cube $\left[-\frac{\mathrm{diam}(K)}{2}, \frac{\mathrm{diam}(K)}{2}\right]^d$. 
Substituting $h(x) = e^{-x\cdot Ax+b\cdot x}p(x)$ in equation \eqref{eq:0_funct} gives the identity 
\begin{equation}
    \begin{split}
        0 & = \int_K e^{-(x-\lambda)\cdot A(x-\lambda)+ b\cdot (x-\lambda)}p(x-\lambda) d\mu(x) \\
         & = \int_{K-k_0} e^{-(y+k_0-\lambda)\cdot A(y+k_0-\lambda) +b\cdot (y+k_0-\lambda)}p(y+k_0-\lambda) \, d\mu(y+k_0) \\
    &= e^{-(k_0-\lambda)\cdot A(k_0-\lambda)+b\cdot (k_0-\lambda)} \\
    & \ \ \ \times \int_{K_0} e^{-y\cdot A y+b\cdot y - (A+A^T)k_0\cdot y} e^{y\cdot (A+A^T)\lambda}p(y+k_0-\lambda) \, d\mu(y+k_0).
    \end{split}
\end{equation}
We divide by $e^{-(k_0-\lambda)\cdot A(k_0-\lambda)+b\cdot (k_0-\lambda)}$ and substitute the complex regular Borel measure $d\nu(y) = e^{-y\cdot A y +b\cdot y - y\cdot (A+A^T)k_0} d\mu(y+k_0)$ on $K_0$ to obtain 
\begin{equation}\label{eq:zzz}
    0  = \int_{K_0} e^{y\cdot (A+A^T)\lambda}p(y+k_0-\lambda) \, d\nu(y).
\end{equation}
Now define
$$
H : \C \to \C, \quad H(z) \coloneqq \int_{K_0} e^{y\cdot (A+A^T)z}p(y+k_0-z) \, d\nu(y).
$$
Then the relation \eqref{eq:zzz} is equivalent to
$$
H(\lambda) = 0, \quad \lambda \in \tfrac{1}{c}\Z^d.
$$
It follows from known results on holomorphic integral transforms \cite[Chap.~XII]{lang}, that $H$ defines an entire function. Denoting by $|\nu|$ the total variation measure of $\nu$ and estimating
\begin{equation}
    \begin{split}
        |H(z)|
        & \leq \int_{K_0} |e^{y\cdot (A+A^T)z}p(y+k_0-z) | \, d|\nu|(y) \\ 
    & \leq C \int_{K_0} e^{\| y \|_\infty \cdot \| (A+A^T)z\|_1} e^{\alpha' \| y+k_0-z \|_1} \, d|\nu|(y)\\
    & \leq C \int_{K_0} e^{ \frac{\mathrm{diam}(K)}{2}  \beta' \| z\|_1} e^{\alpha' \| z \|_1} e^{\alpha' \| y+k_0 \|_1} \, d|\nu|(y)\\
    &  = C e^{(\mathrm{diam}(K) \beta + 2 \alpha) \| z \|_1} \int_{K_0} e^{\alpha' \| y+k_0 \|_1} \, d|\nu|(y)
    \end{split}
\end{equation}
shows that $H \in \mathbf{E}_\sigma(\Cd)$ with $\sigma = \mathrm{diam}(K) \beta + 2 \alpha$.

Since $H$ vanishes on $\Lambda$, it follows from Carlson's theorem that $H$ vanishes identically. By construction of $H$, this implies that the relation \eqref{eq:0_funct} holds for $\tfrac{1}{c}\Z^d$ replaced by $\Rd$. By Lemma \ref{lma:l1}, $\mathcal{T}(h,\Rd)$ is complete in $C(K)$. Hence $\mu=0$ and the statement is proved.
\end{proof}

\begin{remark}[Holomorphic phase retrieval]
    The proof of Proposition \ref{prop:p_ab} leverages techniques from complex analysis, particularly Carlson's theorem. It is noteworthy that there exists a literature on the phase retrieval problem within a purely complex-analytic framework (see, for instance, \cite{perez_2021,Chalendar_Partington_2024,mc2004phase,liehr2023arithmetic}). We previously encountered a specific case of this problem when discussing phase retrieval for the Paley-Wiener space (see Remark \ref{rem:pw}), as every band-limited function extends to an entire function.
\end{remark}

\bibliographystyle{acm}
\bibliography{bibfile}

\end{document}